\documentclass[a4paper,10pt,reqno]{amsart}
\usepackage[a4paper,tmargin=40mm,bmargin=35mm,hmargin=30mm]{geometry}
\usepackage[T1]{fontenc}
\usepackage{lmodern}
\usepackage{amsmath}
\usepackage{amssymb}
\usepackage{amsthm}
\usepackage{stmaryrd}
\usepackage{tikz}
\usepackage[all]{xy}
\usepackage{booktabs}
\usepackage{multirow}
\usepackage{multicol}
\usepackage{here}
\usepackage{aliascnt}
\usepackage{hyperref}
\usepackage[noabbrev]{cleveref}

\newcommand{\NewTheorem}[2]{
	\newaliascnt{#1}{TheoremEnvironment}
	\newtheorem{#1}[#1]{#1}
	\aliascntresetthe{#1}
	\crefname{#1}{#1}{#2}
	\Crefname{#1}{#1}{#2}
}

\theoremstyle{definition}

\NewTheorem{Definition}{Definitions}
\NewTheorem{Remark}{Remarks}
\NewTheorem{Axiom}{Axioms}
\NewTheorem{Example}{Examples}
\NewTheorem{Observation}{Observations}
\NewTheorem{Convention}{Conventions}
\NewTheorem{Notation}{Notations}
\NewTheorem{Setting}{Settings}
\NewTheorem{Question}{Questions}
\NewTheorem{Answer}{Answers}
\NewTheorem{Conjecture}{Conjectures}
\NewTheorem{Problem}{Problems}
\NewTheorem{Solution}{Solutions}
\NewTheorem{Goal}{Goals}
\NewTheorem{Comment}{Comments}
\NewTheorem{Aim}{Aims}
\NewTheorem{Caution}{Cautions}
\NewTheorem{Exercise}{Exercises}
\NewTheorem{Examples}{Examples}
\theoremstyle{plain}
\NewTheorem{Proposition}{Propositions}
\NewTheorem{Lemma}{Lemmas}
\NewTheorem{Theorem}{Theorems}
\NewTheorem{Corollary}{Corollaries}

\crefname{enumi}{}{}
\Crefname{enumi}{}{}
\creflabelformat{enumi}{(#2#1#3)}
\crefname{enumii}{}{}
\Crefname{enumii}{}{}
\creflabelformat{enumii}{(#2#1#3)}
\crefname{enumiii}{}{}
\Crefname{enumiii}{}{}
\creflabelformat{enumiii}{(#2#1#3)}

\makeatletter
\renewcommand{\p@enumii}{}
\renewcommand{\p@enumiii}{}
\makeatother

\numberwithin{equation}{section}
\crefname{equation}{}{}
\Crefname{equation}{}{}
\creflabelformat{equation}{(#2#1#3)}

\newcommand{\SwapSymbols}[1]{
	\expandafter\let\expandafter\temporarysymbol\csname #1\endcsname
	\expandafter\let\csname #1\expandafter\endcsname\csname var#1\endcsname
	\expandafter\let\csname var#1\endcsname\temporarysymbol
}

\SwapSymbols{epsilon}
\SwapSymbols{phi}
\SwapSymbols{Gamma}
\SwapSymbols{Delta}
\SwapSymbols{Theta}
\SwapSymbols{Lambda}
\SwapSymbols{Xi}
\SwapSymbols{Pi}
\SwapSymbols{Sigma}
\SwapSymbols{Upsilon}
\SwapSymbols{Phi}
\SwapSymbols{Psi}
\SwapSymbols{Omega}

\newcommand{\bbN}{\mathbb{N}}

\newcommand{\cC}{\mathcal{C}}
\newcommand{\cD}{\mathcal{D}}

\newcommand{\cN}{\mathcal{N}}

\newcommand{\cS}{\mathcal{S}}
\newcommand{\cT}{\mathcal{T}}
\newcommand{\cU}{\mathcal{U}}

\newcommand{\To}{\longrightarrow}

\DeclareMathOperator{\Hom}{Hom}

\DeclareMathOperator{\Ext}{Ext}

\DeclareMathOperator{\Ann}{Ann}

\DeclareMathOperator{\Ker}{Ker}
\DeclareMathOperator{\Coker}{Coker}
\let\Im\relax
\DeclareMathOperator{\Im}{Im}

\DeclareMathOperator{\Spec}{Spec}
\DeclareMathOperator{\Max}{Max}

\DeclareMathOperator{\Ass}{Ass}
\DeclareMathOperator{\Supp}{Supp}

\title{Cofiniteness with respect to extension of Serre subcategories}

\author{Negar Alipour and Reza Sazeedeh}
\address{Department of Mathematics, Urmia University, P.O.Box: 165, Urmia, Iran}
\email{rsazeedeh@ipm.ir}
\address{Department of Mathematics, Urmia University, P.O.Box: 165, Urmia, Iran}
\email{negaralipur8707@yahoo.com}

\begin{document}

\begin{abstract}
Let $R$ be a commutative noetherian ring, $\frak a$ be an ideal of $R$, $\cS$ be an arbitrary Serre subcategory of $R$-modules satisfying the condition $C_{\frak a}$ and let $\cN$ be the subcategory of finitely generated $R$-modules. In this paper, we define and study $\cN\cS$-$\frak a$-cofinite modules with respect to the extension subcategory $\cN\cS$ as an generalization of the classical notion, namely $\frak a$-cofinite modules. For the lower dimensions, we show that the classical results of $\frak a$-cofiniteness hold for the new notion.
\end{abstract}

\maketitle
\tableofcontents

\section{Introduction}
Throughout this paper $R$ is a commutative noetherian ring, $\frak a$ is an ideal of $R$, $\cS$ is a Serre subcategory of $R$-modules, $N$ is a finitely generated $R$-module and $M$ is an arbitrary $R$-module.  In this paper, we introduce and study the cofiniteness with respect to $\cS$ and $\frak a$. The $R$-module $M$ is said to be $\cS$-$\frak a$-{\it cofinite} if $\Supp M\subseteq V(\frak a)$ and $\Ext_R^i(R/\frak a, M)\in\cS$ for all integers $i\geq 0$. This notion originally goes back to a special case $\cS=\cN$, the subcategory of finitely generated modules, where $\cN$-$\frak a$-cofinite was known as $\frak a$-cofinite, defined for the first time by Hartshorne [H], giving a negative answer to a question of [G, Expos XIII, Conjecture 1.1].  

Our main of this paper is to study the cofiniteness with respect to the extension subcategory $\cN\cS$. The $\cN\cS$-$\frak a$-cofinite modules are the generalization of classical cofinite modules. To be more precise, if $\cS=0$, they are $\frak a$-cofinite modules studied by numerous authors [H, Ma, MV, M1, M2, M3]. When $\cS$ is the subcategory of artinian modules, they are $\frak a$-cominimax modules studied in [Z, BN] and when $\cS$ is the subcategory of all modules of finite support, they are $\frak a$-weakly cofinite modules studied in [DM]. We say that $\cS$ satisfies the condition $C_{\frak a}$ if for every $R$-module $M$, the following implication holds.
\begin{center}
$C_{\frak a}$: If $\Gamma_{\frak a}(M)=M$ and $(0:_M{\frak a})$ is in $\cS$, then
$M$ is in $\cS$.
\end{center}
  
In this paper we assume that $\cS$ satisfies the condition $C_{\frak a}$. In Section 2, we first show if $M$ is an $\cS$-$\frak a$-cofinite $R$-module and $N$ is of dimension $d$, then $\Ext_R^i(N,M)\in\cS$ for each $i\geq 0$ (c.f. \cref{t2}). For an $R$-module $M$, Max$M$ denotes the set of maximal ideals contained in $\Supp_RM$. One of the main results of this section is the following theorem.

\begin{Theorem}
Let $M$ be an $\cN\cS$-$\frak a$-cofinite $R$-module with $\dim M\leq 1$ and let ${\rm Max}M\subseteq \Supp\cS$ (e.g. if $R$ is a local ring). Then $\Ext_R^i(N,M)$ is $\cN\cS$-$\frak a$ cofinite for each $i\geq 0$.
\end{Theorem}
For any non-negative integer $n$, we denote by $\cD_{\leq n}$ the subcategory of all $R$-modules of dimension $\leq n$. It is clear that $\cD_{\leq n}$ is a Serre subcategory of the category of $R$-modules. Let $(R,\frak m)$ be a local ring, let $M$ be a $\cN\cD_{\leq n}$-$\frak a$-cofinite $R$-module with $\dim M\leq 2$ and $\Supp_RM$ be a countable set. Then we show that $\Ext_R^i(N,M)$ is $\cN\cD_{\leq n}$-$\frak a$-cofinite for each $i\geq 0$.

Section 3 is devoted to $\cN\cS$-cofiniteness when $\dim R/\frak a=1$. In this section we assume that ${\rm Max}M\subseteq \Supp\cS$ (e.g. if $R$ is a local ring) and we prove the following theorem which generalizes [M3, Theorem 2.3].
\begin{Theorem}
If $\Supp_RM\subseteq V(\frak a)$, then $M$ is $\cN\cS$-$\frak a$-cofinite if and only if $$\Hom_R(R/\frak a,M), \Ext_R^1(R/\frak a,M)\in\cN\cS.$$ 
\end{Theorem}
 
In \cref{ts2}, we show that the subcategory $\cS(\frak a)=\{M\in R{\rm -Mod}|\hspace{0.1cm}{\rm Max} M\subseteq \Supp\cS\hspace{0.1cm}$and $M$ is $\cN\cS$-$\frak a$-cofinite$\}$ of $R$-modules is abelian. In particular, if $R$ is a local ring, the subcategory of $\cN\cS$-$\frak a$-cofinite modules is abelian.

We end the paper by the following result about $\cN\cS$-$\frak a$-cofiniteness of local cohomology modules which generalizes [NS, Theorem 3.3 and Proposition 3.4]. We have the following theorem. 

\begin{Theorem}
Let $n$ be a non-negative integer. Then $\Ext_R^i(R/\frak a,M)\in\cN\cS$ for all $0\leq i\leq n+1$ if and only if $H_{\frak a}^i(M)$ is $\cN\cS$-$\frak a$-cofinite for all $0\leq i\leq n$ and $\Hom_R(R/\frak a,H_{\frak a}^{n+1}(M))\in\cN\cS$.
\end{Theorem}
For the basic facts about local cohomology, we refer the reader to the textbook by Brodmann and Sharp [BS].

\medskip
\section{Extension of subcategories and cofiniteness}

We denote by $R$-Mod, the category of all $R$-modules. A full subcategory $\cS$ of $R$-Mod is called {\it Serre} if it is
closed under taking submodules, quotients and extensions. Throughout this section $\cS$ is a Serre subcategory of $R$-Mod.

\medskip 

\begin{Lemma}\label{gr}
Let $N$ be a finitely generated $R$-module and $M$ be an arbitrary $R$-module such that for a non-negative integer $n$, we have $\Ext_R^i(N,M)\in\cS$ for all $i\leq n$. Then $\Ext_R^i(L,M)\in\cS$ for any finitely generated $R$-module $L$ with $\Supp_RL\subseteq\Supp_RN$ and all $i\leq n$. 
\end{Lemma}
\begin{proof}
By Gruson's Theorem [V, Theorem 4.1], $L$ admits a finite filtration $$0=L_0\subset L_1\subset\dots \subset L_t=L$$ such that each factor $L_i/L_{i-1}$ is the homomorphic image of a direct sum of finitely many copies of $N$. Using an induction on $t$, we may assume that $t=1$; and hence there is an exact sequence $0\To K\To N^s\To L\To 0$ of $R$-modules. We observe that $\Supp_RK\subseteq \Supp_RN$ and so applying $\Hom_R(-,M)$ and using an induction on $n$, the result follows.
\end{proof}

Let $\frak a$ be an ideal of $R$ and let $\cS$ be a Serre subcategory of $R$-modules. An $R$-module $M$ is said to be $\cS$-$\frak a$-{\it cofinite} if $\Supp M\subseteq V(\frak a)$ and $\Ext_R^i(R/\frak a, M)\in\cS$ for all $i\geq 0$.

\medskip

\begin{Lemma}\label{coff}
Let $x\in\frak a$ and $\Supp_RM\subseteq V(\frak a)$. If $(0:_Mx), M/xM$ are both $\cS$-$\frak a$-cofinite, then so is $M$.
\end{Lemma}
\begin{proof}
Considering $f=x1_M$ and $T^i=\Ext_R^i(R/\frak a,-)$, we have $T^i(f)=\Ext_R^i(R/\frak a,f)=0$ for all $i\geq 0$. We observe that $T^i\Ker f, T^i\Coker f\in\cS$ for all $i\geq 0$. Consequently [M2, Corollary 3.2] implies that $\Ext_R^i(R/\frak a,M)\in\cS$ for all $i\geq 0$.  
\end{proof}
\medskip

\begin{Lemma}\label{me}
 Let $\cS$ be a Serre subcategory of $R$-modules and let $M$ be an $\cS$-$\frak a$-cofinite $R$-module. Then for each $R$-module $N$ of finite length, $\Ext_R^i(N,M)\in\cS$ for each $i\geq 0$.
\end{Lemma}
\begin{proof}
Since $N$ has finite length, there exists a finite filtration $0=N_n\subset N_{n-1}\subset \dots \subset N_1\subset N_0=N$ of submodule of $N$ such that $N_i/N_{i+1}\cong R/\frak m_i$ is simple for $0\leq i\leq n-1$. It suffices to show that $\Ext_R^j(R/\frak m_i,M)\in\cS$ for all $j\geq 0$ and $0\leq i\leq n-1$ and hence we may assume that $N=R/\frak m$ for some maximal ideal $\frak m$ of $R$. If $\Ext_R^i(R/\frak m,M)=0$ for all $i\geq 0$, there is nothing to prove; otherwise, we have $\frak m\in\Supp M\subseteq V(\frak a)$. Then it follows from \cref{gr} that $\Ext_R^i(R/\frak m,M)\in
\cS$ for all $i\geq 0$. 
\end{proof}

\medskip

Given an $R$-module $M$, the subcategory $\cS$ is said to satisfy {\it the condition $C_{\frak a}$
on} $M$ if the following implication holds:
\begin{center}
If $\Gamma_{\frak a}(M)=M$ and $(0:_M{\frak a})$ is in $\cS$, then
$M$ is in $\cS$.
\end{center}
 We say that $\cS$ satisfies the condition $C_{\frak a}$ if $\cS$ satisfy the condition $C_{\frak a}$
on every $R$-module.

In the rest of this section, we may assume that $\frak a$ is an ideal and $\cS$ satisfies the condition $C_{\frak a}$ and we assume that $N$ is a finitely generated $R$-module.

\begin{Theorem}\label{t2}
Let $M$ be an $\cS$-$\frak a$-cofinite $R$-module and let $N$ be of dimension $d$. Then $\Ext_R^i(N,M)\in\cS$ for each $i\geq 0$.
\end{Theorem}
\begin{proof}
We proceed by induction on $d$. If $d=0$, then the result follows by \cref{me} and so we assume that $d>0$. As $\Supp_R\Gamma_{\frak a}(N)\subseteq V(\frak a)$, the assumption and \cref{gr} imply that $\Ext_R^i(\Gamma_{\frak a}(N),M)\in\cS$ for all $i\geq 0$. Thus applying the functor $\Hom_R(-,N)$ to the exact sequence $$0\To\Gamma_{\frak a}(N)\To N\To N/\Gamma_{\frak a}(N)\To 0$$ we may assume that $\Gamma_{\frak a}(N)=0$. Then $\frak a$ contains a non-zero divsior $x$ of $N$ so that there exists an exact sequence of $R$-modules $0\To N\stackrel{x.}\To N\To N/xN\To 0$ such that $\dim N/xN\leq d-1$. Application of $\Hom_R(-,M)$ to the above exact sequence, for each $i\geq 0$, we have an exact sequence $\Ext_R^i(N/xN,M)\To (0:_{\Ext_R^i(N,M)}x)\To 0$. The induction hypothesis implies that $\Ext_R^i(N/xN,M)\in\cS$ and so $(0:_{\Ext_R^i(N,M)}x)\in\cS$ for all $i\geq 0$. Thus $(0:_{\Ext_R^i(N,M)}\frak a)\in\cS$ and since $\cS$ satisfies the condition $C_{\frak a}$, $\Ext_R^i(N,M)\in\cS$ for all $i\geq 0$.    
\end{proof}

\medskip

\begin{Corollary}
Let $R$ be a local ring and let $M$ be an $\cS$-$\frak a$-cofinite $R$-module. Then $\Ext^i_R(N,M)\in\cS$ for each $i\geq 0$.
\end{Corollary}
\begin{proof}
Since $R$ is local, every finitely generated $R$-module has finite Krull dimension; and hence the result follows by \cref{t2}. 
\end{proof}

\medskip
For a Serre subcategory $\cS$ of $R$-modules, the support of $\cS$ is denoted by $\Supp\cS$ which is $\Supp\cS=\bigcup_{M\in\cS}\Supp_RM=\{\frak p\in\Spec R|\hspace{0.1cm} R/\frak p\in\cS\}.$ The full subcategory of finitely generated $R$-modules is denoted by $\cN$. We denote by $\cN\cS$, the extension subcategory of $\cN$ and $\cS$ which is: $$\cN\cS=\{M\in\cC|\hspace{0.1cm} {\rm there \hspace{0.1cm} exists \hspace{0.1cm} an\hspace{0.1cm}exact\hspace{0.1cm}
 sequence}\hspace{0.1cm}
0\To N\To M\To S\To 0\hspace{0.1cm} {\rm with}\hspace{0.1cm} N\in\cN
\hspace{0.1cm}{\rm and}\hspace{0.1cm} S\in\cS\}.$$
If $\cS$ is a Serre subcategory of $R$-Mod, then by virtue of [Y, Corollary 3.3], $\cN\cS$ is Serre. 

\medskip

\begin{Corollary}
Let $R/\frak a\in\cS$, let $M$ be an $\cN\cS$-$\frak a$-cofinite $R$-module and let $N$ be of dimension $d$. Then $\Ext_R^i(N,M)\in\cN\cS$ for each $i\geq 0$. 
\end{Corollary}
\begin{proof}
Since $\cS$ satisfies the condition $C_{\frak a}$, it follows from [AMS, Theorem 3.8] that $\cN\cS$ satisfies the condition $C_{\frak a}$. Now, the result follows from \cref{t2}.
\end{proof}

For any ideal $\frak a$ of $R$, {\it arithmetic rank of} $R$, denoted by ara$\frak a$, is the least non-negative integer of elements of $R$ required to generate an ideal which has the same radical as $\frak a$. Thus 
$${\rm ara}\frak a=\min\{n\in\bbN_0|\hspace{0.1cm}\exists\hspace{0.1cm} a_1,\dots,a_n\in R \hspace{0.1cm}\textnormal{with}\hspace{0.1cm} \sqrt{(a_1,\dots,a_n)}=\sqrt{\frak a}\}.$$

  For every $R$-module $M$, ara$_M\frak a$ is the arithmetic rank of the ideal $\frak a+\Ann_RM/\Ann_RM$ of the ring $R/\Ann_RM$. We denote by Max$M$ the set of maximal ideals in $\Supp_RM$.

\medskip

\begin{Theorem}\label{t4}
Let $M$ be an $\cN\cS$-$\frak a$-cofinite $R$-module with $\dim M\leq 1$ and ${\rm Max}M\subseteq \Supp\cS$ (e.g. if $R$ is a local ring). Then $\Ext_R^i(N,M)$ is $\cN\cS$-$\frak a$ cofinite for each $i\geq 0$.
\end{Theorem}
\begin{proof}
We proceed by induction on $n={\rm ara}_N\frak a={\rm ara}(\frak a+{\rm Ann}_RN/{\rm Ann}_RN)$. If $n=0$, then there exists some positive integer $t$ such that $N=(0:_N\frak a^t)$ and so the result follows from \cref{gr}. As $\Ann_RN\subseteq\Ann_R N/\Gamma_{\frak a}(N)$, we have ${\rm ara}_{N/\Gamma_{\frak a}(N)}\frak a\leq {\rm ara}_{N}{\frak a}$ and so considering the exact sequence $$0\To \Gamma_{\frak a}(N)\To N\To N/\Gamma_{\frak a}(N)\To 0$$ and \cref{gr}, we may assume that $\Gamma_{\frak a}(N)=0$. If $\Phi=\{\frak p\in\Ass_RM\cap \Supp \cS|\hspace{0.1cm}\dim R/\frak p=1\}$, then using [B, Ch. IV, Sec.1.2, Proposition 4], there exists a submodule $K$ of $M$ such that $\Ass_RK=\Phi$ and $\Ass_RM/K=\Ass_RM\setminus \Phi$. Since $M$ be is $\cN\cS$-$\frak a$-cofinite, $\Hom_R(R/\frak a,K)\in\cN\cS$ and so there is an exact sequence of $R$-modules $$0\To F\To\Hom_R(R/\frak a,K)\To S\To 0$$ such that $F$ is finitely generated and $S\in\cS$. Every $\frak q\in\Supp F$ contains a prime ideal $\frak p\in\Ass K$ and hence there is an epimorphism $R/\frak p\To R/\frak q\To 0$. The fact that $R/\frak p\in\cS$ implies that $R/\frak q\in\cS$. Since $F$ is noetherian, there is a finite filtration of submodules of $F$ $$0=F_m\subseteq F_{m-1}\subseteq\dots F_1\subseteq F_0=F$$ and prime ideals $\frak p_i\in\Supp F, 0\leq i\leq m-1$ such that $N_i/N_{i+1}\cong R/\frak p_i\in\cS$. This forces that $F\in\cS$; and hence $\Hom_R(R/\frak a,K)\in\cS$. Since $\cS$ satisfies the condition $C_{\frak a}$, we deduce that $K\in\cS$. Thus for every finitely generated $R$-module $L$, the module $\Ext_R^i(L,K)\in\cS$ for all $i\geq 0$. Therefore, replacing $M$ by $M/K$ we may assume that every $\frak p\in\Ass_RM$ with $\dim R/\frak p=1$ is not in Supp$\cS$. For a non-negative integer $t$, let $\cT_t=\bigcup_{i=0}^t\Supp\Ext_R^i(N,M)$ and $\cT=\{\frak p\in\cT_t|\hspace{0.1cm} \dim R/\frak p=1\}$. We notice that
 $\{\frak p\in\Ass_RM|\hspace{0.1cm}\dim R/\frak p=1\}$ is a finite set and $\cT\subseteq\{\frak p\in\Ass_RM|\hspace{0.1cm}\dim R/\frak p=1\}$ and hence $\cT$ is a finite set. The assumption implies that $\Hom_R(R/\frak a,M)\in\cN\cS$ so that there exists an exact sequence $0\To F\To \Hom_R(R/\frak a,M)\To S\To 0$ of $R$-modules such that $F$ is finitely generated and $S\in\cS$. For every $\frak p\in\cT$, since $\frak p\notin\Supp\cS$, localizing at $\frak p$, the $R_{\frak p}$-module $\Hom_R(R/\frak a,M)_{\frak p}\cong F_{\frak p}$ has finite length so that $M_{\frak p}$ is an artinian  and $\frak a$-cofinite by [M1, Theorem 1.6]. It therefore follows from [M1, Corollary 1.7] that $\Ext_R^i(N,M)_{\frak p}$ is artinian and $\frak aR_{\frak p}$-cofinite for all $i\geq 0$. Let $\cT=\{\frak p_1, \dots, \frak p_l\}$. By [BN, Lemma 2.5], for all $0\leq i\leq k$ and all $1\leq j\leq n$, we have 
$$V(\frak aR_{\frak p_j})\cap {\rm Att}_{R_{\frak p_j}}(\Ext_R^i(N,M))_{\frak p_j}\subseteq V(\frak p_jR_{\frak p_i}).$$If we set $\cU=\bigcup_{i=0}^k\bigcup_{j=1}^l\{\frak q\in\Spec R|\hspace{0.1cm}\frak qR_{\frak p_j}\in{\rm Att}_{R_{\frak p_j}}(\Ext_R^i(N,M))_{\frak p_j}\}$ for all $0\leq i\leq k$ and all $1\leq j\leq l$, then $\cU\cap V(\frak a)\subseteq \cT$. For each $i\geq 0$, we have $\Ann_RN\subseteq \Ann\Ext_R^i(N,M)$; and hence for every $\frak q\in\cU$, we have $(\Ann_RN)R_{\frak p_j}\subseteq \frak qR_{\frak p_j}$ where $\frak qR_{\frak p_j}\in{\rm Att}_{R_{\frak p_j}}(\Ext_R^i(N,M))$ for some $0\leq i \leq k $ and $1\leq j\leq l$. This implies $\Ann_RN\subseteq \frak q$ so that $\cU\subseteq \Supp N$. Since ${\rm ara}_N\frak a=n$, there exists $a_1,\dots,a_n\in R$ such that $\sqrt{\frak a+\Ann_RN}=\sqrt{(a_1,\dots,a_n)+\Ann_RN}$. Since $\frak a\nsubseteq (\bigcup_{\frak q\in\cU\setminus V(\frak a)}\frak q)\bigcup (\bigcup_{\frak p\in\Ass N}\frak p)$, we deduce that $(y_1,\dots,y_n)\nsubseteq (\bigcup_{\frak q\in\cU\setminus V(\frak a)}\frak q)\bigcup (\bigcup_{\frak p\in\Ass N}\frak p)$ and so using [M, Exercise 16.8], there exists $b\in(y_2,\dots,y_n)$ such that $x=y_1+b\notin (\bigcup_{\frak q\in\cU\setminus V(\frak a)}\frak q)\bigcup (\bigcup_{\frak p\in\Ass N}\frak p)$. It is clear that $(y_1,\dots,y_n)=(x,y_2,\dots,y_n)$ and so $(y_1,\dots,y_n)+\Ann_RN/xN=(y_2,\dots,y_n)+\Ann_RN/xN)$. Thus ${\rm ara}_{N/xN}\frak a\leq n-1$ and there is an exact sequence of $R$-modules $0\To N\stackrel{x.}\To N\To N/xN\To 0$ which induces the following exact sequence of $R$-modules $$\Ext_R^{i}(N/xN,M)\To\Ext_R^i(N,M)\stackrel{x.}\To\Ext_R^i(N,M)\To \Ext_R^{i+1}(N/xN,M).$$  If we consider $D_i=\Ext_R^i(N/xN,M)$ and $L_i=\Ext^i_R(N,M)/x\Ext_R^i(N,M)$, using the induction hypothesis, $D_i$ is $\cN\cS$-$\frak a$-cofinite for all $i\geq 0$. On the other hand, it follows from [BN, Lemma 2.4] that $(L_i)_{\frak p_j}$ has finite length; and hence there exists a finitely generated submodule $L_{ij}$ of $L_i$ such that $(L_i)_{\frak p_j}={L_{ij}}_{\frak p_j}$ for each $0\leq i\leq t$ and $1\leq j\leq l$. For each $0\leq i\leq t$, let $L'_i=L_{i1}+\dots+L_{il}$. Then $L'_i$ is a finitely generated submodule of $L$ and so the previous argument and the assumption on $M$ imply that $\Supp_RL_i/L'_i\subseteq \cT_t\setminus\cT\subseteq\Max R\bigcap\Supp\cS$. We prove that $L_i\in\cN\cS$ for all  $0\leq i\leq t$. Since $D_{i+1}/L'_i$ is $\cN\cS$-$\frak a$-cofinite and $L_i/L'_i$ is a submodule of $D_{i+1}/L'_i$, the module $\Hom_R(R/\frak a, L_i/L'_i)\in\cN\cS$. Then there exists an exact sequence of $R$-modules 
$$0\To F\To \Hom_R(R/\frak a, L_i/L'_i)\To S\To o$$ such that $F$ is finitely generated and $\cS\in\cS$. Since $\Supp_R \Hom_R(R/\frak a, L_i/L'_i)\subseteq{\rm Max}R\cap\Supp\cS$, the module $F$ has finite length and $F\in\cS$ so that $\Hom_R(R/\frak a,L_i/L'_i)\in\cS$. Since $\cS$ satisfies the condition $C_{\frak a}$, we deduce that $L_i/L'_i\in\cS$. This implies that $L_i\in\cN\cS$ for all $0\leq i\leq t$ and the exact sequence $$0\To L_i\To D_{i+1}\To (0:_{\Ext_R^{i+1}(N,M)}x)\To 0$$ implies that $(0:_{\Ext_R^{i}(N,M)}x)$ is $\cN\cS$-$\frak a$-cofinite for all $1\leq i\leq t$. Moreover, $(0:_{\Hom_R(N,M)}x)\cong\Hom_R(N/xN,M)$ is $\cN\cS$-$\frak a$-cofinite by the induction hypothesis. It now follows from \cref{coff} that $\Ext_R^i(N,M)$ is $\cN\cS$-$\frak a$-cofinite for all $0\leq i\leq t$. Since $t$ is arbitrary, we deduce that $\Ext_R^i(N,M)$ is $\cN\cS$-$\frak a$-cofinite for all $i\geq 0$.
 \end{proof}

For any non-negative integer $n$, we denote by $\cD_{\leq n}$ the subcategory of all $R$-modules of dimension $\leq n$. It is clear that $\cD_{\leq n}$ is a Serre subcategory of the category of $R$-modules. 

\medskip

\begin{Corollary}\label{c5}
Let $n$ be a non-negative integer and let $M$ be a $\cN\cD_{\leq n}$-$\frak a$-cofinite $R$-module with $\dim M\leq 1$. Then $\Ext_R^i(N,M)$ is $\cN\cD_{\leq n}$-$\frak a$-cofinite for each $i\geq 0$.
\end{Corollary}
\begin{proof}
 It is clear that $\cD_{\leq n}$ satisfies the condition $C_{\frak a}$ for all ideal $\frak a$ of $R$ and so the result follows by \cref{t4}. 
\end{proof}
\medskip
\begin{Corollary}
Let $(R,\frak m)$ be a local ring, let $M$ be a $\cN\cD_{\leq n}$-$\frak a$-cofinite $R$-module with $\dim M\leq 2$ and a non-negative integer $n$, and let $\Supp_{\hat{R}}(M\otimes_R\hat{R})$ be a countable set. Then $\Ext_R^i(N,M)$ is $\cN\cD_{\leq n}$-$\frak a$-cofinite for each $i\geq 0$.
\end{Corollary}
\begin{proof}
In view of \cref{t4}, it suffices to consider that $\dim M=2$. There exists a prime ideal $\frak p\in\Ass M$ such that $\dim R/\frak p=\dim \hat{R}/\frak p\hat{R}=2$ where $\hat{R}$ is the completion of $R$ with respect to $\frak m$-adic-topology. Since $R/\frak p$ is a submodule of $M, \hat{R}/\frak p\hat{R}$ is a submodule of $M\otimes_R\hat{R}$ so that $\dim _{\hat{R}}(M\otimes_R\hat{R})\geq 2$. If $\dim _{\hat{R}}(M\otimes_R\hat{R})=t$ for some $t$, there exists $\frak P\in\Ass _{\hat{R}}(M\otimes_R\hat{R})$ such that $\dim \hat{R}/\frak P=t$ and $\frak P=\Ann_{\hat{R}}(x)$ where $x\in M\otimes_R\hat{R}$. Then there exists a finitely generated submodule $K$ of $M$ such that $\frak P\in\Ass_{\hat{R}}(K\otimes _R\hat{R})$. But $t=\dim _{\hat{R}}(K\otimes _R\hat{R})=\dim_RK\leq 2$ and hence $\dim_{\hat{R}}(M\otimes_R\hat{R})=2$. Since $M$ is $\cN\cD_{\leq n}$-$\frak a$-cofinite, for each $i\geq 0$, there exists an exact sequence of $R$-modules $0\To K\To \Ext_R^i(R/\frak a,M)\To D\To 0$ such that $K$ is finitely generated and $\dim D\leq n$. A similar argument mentioned above, implies that $\dim_{\hat{R}}(D\otimes_R\hat{R})\leq n$ and so $M\otimes_R\hat{R}$ is $\hat{\cN}\hat{\cD}$-$\frak a\hat{R}$-cofinite where $\hat{\cN}$ denotes the subctegory of finitely generated $\hat{R}$-modules and $\hat{\cD}_{\leq n}$ denotes the subcategory of all $R$-modules of dimension $\leq n$. For each $i\geq 0$, if $\Ext_{\hat{R}}^i(N\otimes_R\hat{R},M\otimes_R\hat{R})\cong \Ext_R^i(N,M)\otimes_R\hat{R}$ is a $\hat{\cN}\hat{\cD}_{\leq n}$-$\frak a\hat{R}$-cofinite module, then for each $j\geq 0$, there exists an exact sequence of $\hat{R}$-modules
$$0\To X\To \Ext_R^j(R/\frak a,\Ext_R^i(N,M))\otimes_R\hat{R}\To Y\To 0$$ such that $X$ is finitely generated and $\dim Y\leq n$. It is clear that there exits a finitely generated $R$-submodule $N$ of $\Ext_R^j(R/\frak a,\Ext_R^i(N,M))$ such that $X= N\otimes_R\hat{R}$ and hence $Y\cong (\Ext_R^j(R/\frak a,\Ext_R^i(N,M))/N)\otimes_R\hat{R}$ so that $\dim\Ext_R^j(R/\frak a,\Ext_R^i(N,M))/N\leq n$ by a similar argument mentioned in the beginning of the proof . This implies that $\Ext_R^i(N,M)$ is $\cN\cS$-$\frak a$-cofinite for all $i\geq 0$. On the other hand, by virtue of [Ma, Lemma 2.1], we have $$\Supp_RM=\bigcup_{K\leq M}\Ass_RM/K\subseteq\{\frak p\cap R|\hspace{0.1cm} \frak p\in\Ass_{\hat{R}}(M\otimes_R\hat{R}/K\otimes_R\hat{R}\}\subseteq\{\frak p\cap R|\hspace{0.1cm}\frak p\in
\Supp_{\hat{R}}(M\otimes_R\hat{R})\}$$ 
which implies that $\Supp_RM$ is a countable set. Then without loss of generality we may assume that $R$ is complete. If we consider $\cT=\{\frak p\in\Supp_RM|\hspace{0.1cm} \dim R/\frak p=1\}$, then it follows from [MV, Lemma 3.2] that $\frak m\nsubseteq\bigcup_{\frak p\in\cT}\frak p$. Letting $S=R\setminus\bigcup_{\frak p\in\cT}\frak p$, it is clear that $\dim_{S^{-1}R} S^{-1}M\leq 1$ and $S^{-1}M$ is an $\cN'\cD'_{\leq n-1}$-$S^{-1}\frak a$-cofinite $S^{-1}R$-module where $\cN'$ is the subcategory of finitely generated $S^{-1}R$-modules and $\cD'_{n-1}$ is the subcategory of all $S^{-1}R$-modules of dimension $\leq n-1$. Then, in view of \cref{c5}, for any finitely generated $R$-module $N$, the $S^{-1}R$-module $\Ext_{S^{-1}}^i(S^{-1}N,S^{-1}M)$ is $\cN'\cD'_{\leq n-1}$-$S^{-1}\frak a$-cofinite for each $i\geq 0$. Thus for each $i\geq 0$ and each $j\geq 0$, there is an exact sequence of $S^{-1}R$-modules $$0\To N'\To S^{-1}\Ext_{R}^j(R/\frak a,\Ext_{R}^i(N,M))\To D'\To 0$$ such that $N'$ is finitely generated and $D'\in\cD'_{n-1}$. Whence, there is a finitely generated submodule $N$ of $\Ext_{R}^j(R/\frak a,\Ext_{R}^i(N,M))$ such that $S^{-1}N=N'$ and $D'=S^{-1}D$ where $D=\Ext_{R}^j(R/\frak a,\Ext_{R}^i(N,M))/N\in\cD_n.$. Consequently, $\Ext_{R}^j(R/\frak a,\Ext_{R}^i(N,M))\in\cN\cD_n.$ 
\end{proof}

\medskip

 
\section{Cofiniteness with respect an ideal of dimension one} 
Throughout this section $\frak a$ is an ideal of $R$ with $\dim R/\frak a=1$ and $\cS$ is a Serre subcategory of $R$-modules satisfying the condition $C_{\frak a}$.
\begin{Lemma}\label{l1}
Let $M$ be an $R$-module such that $\Supp M\subseteq V(\frak a)$ and $\Ass_RM\cap\Supp\cS\subseteq{\rm Max}R$. If $\Hom_R(R/\frak a,M)\in\cN\cS$, then there is a finitely generated submodule $N$ of $M$ and an element $x\in\frak a$ such that $\Supp_R (M/(xM+N))\subseteq\Max R$. 
\end{Lemma}
\begin{proof}
By the assumption, there exists an exact sequence of $R$-modules $$0\To N\To \Hom_R(R/\frak a,M)\To S\To 0$$ such that $N$ is finitely generated and $S\in\cS$. We observe that $\Supp_RS\subseteq{\rm Max}R$ because if $\frak q\in\Supp\cS$ is a non-maximal ideal of $R$, then $\dim R/\frak q=1$ so that $\frak q\in\Ass_RM$ which is a contradiction by the assumption. Since $\dim R/\frak a=1$, there exists finitely many prime ideals $\frak p_1,\dots,\frak p_n$ containing $\frak a$. Considering $T=R\setminus\cup_{i=1}^n\frak p_i$, we have $T^{-1}N=(0:_{T^{-1}M}T^{-1}\frak a)$ is a finitely generated $T^{-1}R$-module. Using a similar proof of [M3, Proposition 2.2], there exists an element $x\in\frak a$ and a finitely generated submodule $N$ of $M$ such that
$\Supp_R (M/(xM+N))\subseteq\Max R$.
\end{proof}

The following theorem generalizes [M3, Theorem 2.3].
\medskip

\begin{Theorem}\label{ts1}
Let $M$ be an $R$-module such that $\Supp_RM\subseteq V(\frak a)$ and ${\rm Max}M\subseteq \Supp\cS$ (e.g. if $R$ is a local ring). Then $M$ is $\cN\cS$-$\frak a$-cofinite if and only if $\Hom_R(R/\frak a,M), \Ext_R^1(R/\frak a,M)\in\cN\cS$. 
\end{Theorem}
\begin{proof}
A part of the proof is similar to the proof of [M3, Proposition 2.3]. If the theorem does not hold, there is an $R$-module $M$ whose annihilator is maximal among those ideals, which occurs as annihilator of $R$-modules satisfying the hypothesis, but are not $\cN\cS$-$\frak a$-cofinite.  
Let $\Phi=\{\frak p\in\Ass_RM|\hspace{0.1cm} \dim R/\frak p=1\}\cap\Supp\cS$. In view of [B, Chap. IV. Sec 1.2, Proposition 4], there exists a submodule $K$ of $M$ such that $\Ass_RK=\Phi$ and $\Ass_RM/K=\Ass_RM\setminus \Phi$. We observe by the assumption that $\Hom_R(R/\frak a,K)\in\cN\cS$ and so there exists an exact sequence of $R$-modules $$0\To N\To \Hom_R(R/\frak a,K)\To S\To 0$$ such that $N$ is finitely generated and $S\in\cS$. Considering a finite filtration of $N$ and the fact that $\Ass_RN\subseteq\Supp\cS$, we deduce that $N\in\cS$ and so $\Hom_R(R/\frak a,K)\in\cS$. Since $\cS$ satisfies the condition $C_{\frak a}$, we have $K\in\cS$. Therefore, replacing $M$ by $M/K$, we may assume that for every $\frak p\in\Ass_RM$ with $\dim R/\frak p=1$, we have $\frak p\notin\Supp\cS$; and hence $\Ass_RM\cap\Supp\cS\subseteq {\rm Max}R$. Since $\Hom_R(R/\frak a,M)\in\cN\cS$, it follows from \cref{l1} that there exists $x\in\frak a$ and a finitely generated submodule $N$ of $M$ such that $\Supp_R (M/(xM+N))\subseteq\Max R$. We observe that $M/N$ satisfies the hypothesis and $M$ is $\cN\cS$-$\frak a$-cofinite if and only if $M/N$ is $\cN\cS$-$\frak a$-cofinite and the inclusion $\Ann_RM\subseteq\Ann_RM/N$ is equal. Then we can replace $M$ by $M/N$ and we may assume that $\Supp_R(M/xM)\subseteq \Max R$. If $xM=0$, we have $\Supp _RM\subseteq \Max R$ and so by the assumption we have $\Supp_RM\subseteq\Supp\cS$. Since $\Hom_R(R/\frak a,M)\in\cN\cS$, there exists an exact sequence of $R$-modules $$0\To N\To \Hom_R(R/\frak a,M)\To \cS\To 0$$ such that $N$ is finitely generated and $S\in\cS$. 
It is clear that $N$ has finite length and the fact that $\Supp_RM\subseteq\Supp\cS$ and the previous argument implies that $N\in\cS$, and hence $\Hom_R(R/\frak a,M)\in\cS$. Since $\cS$ satisfies the condition $C_{\frak a}$, we have $M\in\cS$ so that $M$ is $\cN\cS$-$\frak a$-cofinite which is a contradiction. Then $x\notin\Ann_RM$. Considering the exact sequences $$0\To (0:_Mx)\To M\To xM\To 0;$$  $$0\To xM\To M\To M/xM\To 0,$$ it is clear that $\Hom_R(R/\frak a,(0:_Mx)),\Ext_R^1(R/\frak a,(0:_Mx))\in\cN\cS$ and $\Ann_RM\subsetneq\Ann_R(0:_Mx)$. The maximality implies that $(0:_Mx)$ is $\cN\cS$-$\frak a$-cofinite. The exact sequences imply that $\Hom_R(R/\frak a,M/xM)\in\cN\cS$ and by the above argument and the assumption, we have $\Supp_RM/xM \subseteq{\rm Max}R\cap\Supp\cS$. Using a similar argument mentioned before and the fact that $\cS$ satisfies the condition $C_{\frak a}$, we deduce that $M/xM\in\cS$ so that $M/xM$ is $\cN\cS$-$\frak a$-cofinite. Consequently, \cref{coff} implies that $M$ is $\cN\cS$-$\frak a$-cofinite which is a contradiction.   
\end{proof}

\medskip

\begin{Corollary}
Let $M$ be an $R$-module with $\Supp_RM\subseteq V(\frak a)$ and $\Hom_R(R/\frak a,M),\Ext_R^1(R/\frak a,M)\in\cN\cS$, let ${\rm Max}M\subseteq \Supp\cS$ (e.g. if $R$ is a local ring), and let $N$ be a finitely generated $R$-module. Then $\Ext_R^i(N,M)$ is $\cN\cS$-$\frak a$ cofinite for each $i\geq 0$.
\end{Corollary}
\begin{proof}
Since $\Supp_RM\subseteq V(\frak a)$, we have $\dim M\leq 1$. Now the result is obtained by \cref{t4} and \cref{ts1}.
\end{proof}

\medskip
The following theorem generalizes [M3, Theorem 2.6]. 

\begin{Theorem}\label{ts2}
The subcategory $\cS(\frak a)=\{M\in R{\rm -Mod}|\hspace{0.1cm}{\rm Max} M\subseteq \Supp\cS\hspace{0.1cm}$and $M$ is $\cN\cS$-$\frak a$-cofinite$\}$ of $R$-modules is abelian. In particular, if $R$ is a local ring, the subcategory of $\cN\cS$-$\frak a$-cofinite modules is abelian.
\end{Theorem}
\begin{proof}
Given an $R$-homomorphism $f:M\To N$ in $\cS(\frak a), K=\Ker f, I=\Im f$ and $C=\Coker f$, it is straightforward to show that $\Hom_R(R/\frak a,K),\Ext_R^1(R/\frak a,K)\in\cN\cS$ and hence using \cref{ts1}, the module $K$ is $\cN\cS$-$\frak a$-cofinite. This implies that $I$ and consequently $C$ is $\cN\cS$-$\frak a$-cofinite.
\end{proof}
\medskip
For $\cN\cS$-$\frak a$-cofiniteness of local cohomology modules, we have the following theorem which generalizes 
[NS, Theorem 3.3 and Proposition 3.4].
\begin{Theorem}\label{ts3}
Let $M$ be an $R$-module such that ${\rm Max}M\subseteq\Supp\cS$ (e.g. if $R$ is a local ring) and  let $n$ be a non-negative integer. Then $\Ext_R^i(R/\frak a,M)\in\cN\cS$ for all $i\leq n+1$ if and only if $H_{\frak a}^i(M)$ is $\cN\cS$-$\frak a$-cofinite for all $i\leq n$ and $\Hom_R(R/\frak a,H_{\frak a}^{n+1}(M))\in\cN\cS$.
\end{Theorem}
\begin{proof}
We show bi-implication by induction on $n$. If $n=0$ and $\Ext_R^i(R/\frak a,M)\in\cN\cS$ for $i=0,1$. It is straightforward to see that $\Hom_R(R/\frak a,\Gamma_{\frak a}(M)), \Ext_R^1(R/\frak a,\Gamma_{\frak a}(M))\in\cN\cS$; and hence according to \cref{ts1}, the module $\Gamma_{\frak a}(M)$ is $\cN\cS$-$\frak a$-cofinite.  On the other hand, there exists an exact sequence of $R$-modules $0\To M/\Gamma_{\frak a}(M)\To E\To Q\To 0$ such that $E$ is injective with $\Gamma_{\frak a}(E)=0$. Thus in view of the exact sequence of $R$-modules $$0\To \Gamma_{\frak a}(M)\To M\To M/\Gamma_{\frak a}(M)\To 0\hspace{0.2cm}(\dag)$$ we have the following isomorphims 
$$\Hom_R(R/\frak a,H_{\frak a}^1(M))\cong \Hom_R(R/\frak a,\Gamma_{\frak a}(M))\cong\Hom_R(R/\frak a,Q)\cong \Ext_R^1(R/\frak a,M/\Gamma_{\frak a}(M))\in\cN\cS.$$
 Conversely, $\Hom_R(R/\frak a,M)\cong \Hom_R(R/\frak a,\Gamma_{\frak a}(M))\in\cN\cS$ by the assumption. Furtheremore, by the above isomorphisms, we have $\Ext_R^1(R/\frak a,M/\Gamma_{\frak a}(M))\in\cN\cS$; and hence the exact sequence $(\dag)$ implies that $\Ext_R^1(R/\frak a,M)\in\cN\cS$. 
Assume that $n>0$ and so by the induction step, $\Gamma_{\frak a}(M)$ is $\cN\cS$-$\frak a$-cofinite. Thus the exact sequence $(\dag)$ implies that $\Ext_R^i(R/\frak a,M)\in\cN\cS$ if and only if $\Ext_R^i(R/\frak a,M/\Gamma_{\frak a}(M))\in\cN\cS$  and $H_{\frak a}^i(M)\cong H_{\frak a}^i(M/\Gamma_{\frak a}(M))$ for all $0\leq i\leq n+1$. Then we may assume that $\Gamma_{\frak a}(M)=0$; and hence there is an exact sequence of $R$-modules $$0\To M\To E\To Q\To 0$$ such that $E$ is injective with $\Gamma_{\frak a}(E)=0$. The induction hypothesis implies that $\Ext_R^i(R/\frak a,Q)\in\cN\cS$ for all $0\leq i\leq n$ if and only if $H_{\frak a}^{n-1}(Q)$ is $\cN\cS$-$\frak a$-cofinite if and $\Hom(R/\frak a,H_{\frak a}^n(Q))\in\cN\cS$. Consequently the isomorphisms $\Ext_R^i(R/\frak a,Q)\cong\Ext_R^{i+1}(R/\frak a,M)$ and $H_{\frak a}^i(Q)\cong H_{\frak a}^{i+1}(M)$ for all $i\geq 0$ get the assertion. 
\end{proof}



\end{document}